\documentclass[12pt,a4paper]{article}
\usepackage{amsfonts,amsmath,amssymb,indentfirst}
\newcommand{\er}{\mathbb{R}}

\newtheorem{Teorema}{Theorem}[section]

\newtheorem{Rem}[Teorema]{Remark}

\begin{document}
\title{Ground states for a coupled nonlinear Schr\"odinger system}
\author{Filipe Oliveira\footnote{e-mail address: fso@fct.unl.pt}}
\date{Department of Mathematics, FCT-UNL-CMA\\Universidade Nova de Lisboa\\
Caparica Campus, 2829-516, Portugal}

\maketitle

\begin{abstract}
\noindent
 We study the existence of ground states for the coupled Schr\"odinger system
 \begin{equation}
\label{ellipticabstract}
\left\{
\begin{array}{llll}
-\Delta u+u&=&|u|^{2q-2}u+b|v|^q|u|^{q-2}u\\
-\Delta v+\omega^2v&=&|v|^{2q-2}v+b|u|^q|v|^{q-2}v
\end{array}\right.
\end{equation}
\noindent
in $\er^n$, for $\omega \geq 1$, $b>0$ (the so-called ``attractive case'') and $q>1$ ($q<\frac n{n-2}$ if $n\geq 3$). We improve for  several ranges of $(q,n,\omega)$ the known results concerning the existence of positive ground state solutions to (\ref{ellipticabstract}) with non-trivial components. In particular, we prove that for $1<q<2$ such ground states exist in all dimensions and for all values of $\omega$, which constitutes a drastic change of behaviour with respect to the case $q\geq 2$. Furthermore, for $q>2$ and in the one-dimensional case $n=1$, we improve the results in \cite{MMP}.
\end{abstract}

\section{Introduction}
In this paper we consider the system
\begin{equation}
\label{initial}
\left\{
\begin{array}{llll}
-\Delta u+\lambda_1u&=&|u|^{2q-2}u+b|v|^q|u|^{q-2}u\\
-\Delta v+\lambda_2v&=&|v|^{2q-2}v+b|u|^q|v|^{q-2}v,
\end{array}\right.
\end{equation}
with $u,v\,:\er^n\to\er$ ($n\geq 1$), $q>1$, $b>0$ and $\lambda_1, \lambda_2>0$, which appears in several physical contexts, namely in nonlinear optics (see \cite{Ph1} and the references therein).\\

By rescaling the $x$ variable and/or inverting the roles of $u$ and $v$, it is easy to see that (\ref{initial}) can be reduced, without loss of generality, to the system 
\begin{equation}
\label{elliptic}
\left\{
\begin{array}{llll}
-\Delta u+u&=&|u|^{2q-2}u+b|v|^q|u|^{q-2}u\\
-\Delta v+\omega^2v&=&|v|^{2q-2}v+b|u|^q|v|^{q-2}v,\qquad \omega\geq 1.
\end{array}\right.
\end{equation}
In the last years, this system has been extensively studied by many authors (see for instance \cite{N4}, \cite{N1}, \cite{N2}, \cite{N3}). In particular, in \cite{Colorado} and \cite{Coloradobis} the authors studied the case $q=2$ and $n=2,3$, proving the existence of a constant $\Lambda>0$ depending on $\omega$ such that for $b<\Lambda$ the system (\ref{elliptic}) admits a non-trivial radial solution $(u,v)\neq (0,0)$ (with $u,v>0$ if $b>0$). The authors also showed the existence of another constant $\Lambda'\geq \Lambda$ such that for $b>\Lambda'$ the system possesses a radial ground state solution $W_*=(u_*,v_*)$ ($u_*,v_*>0$), in the sense that $W_*$ minimizes the energy functional
associated to (\ref{elliptic}) among all solutions in $(u,v)\in H^1(\er^n)\times H^1(\er^n)\setminus\{(0,0)\}$. In \cite{IT} Ikoma and Tamaka showed that for $0<b<\min\{\Lambda,1\}$, the solutions found in \cite{Colorado},\cite{Coloradobis} are in fact also least energy solutions. 

\bigskip

\bigskip

In \cite{MMP}, following some of the ideas presented in \cite{Rab}, the authors proved the existence of a radial non-trivial ground state solution $(u^*,v^*)$ ($u^*,v^*\geq 0$) for every $b>0$ and for $(q,n)$ satisfying 
\begin{equation}
 \label{cond}
 1<q<\left\{\begin{array}{llll}
            +\infty&if& n=1,2\\
            \\
            \displaystyle\frac{n}{n-2}&if& n\geq 3.
           \end{array}\right.
\end{equation}
Furthermore, it is shown that for 
\begin{equation}
 \label{cond2}
 b\geq \mathcal{C}_{\omega,n,q}:=\frac 12\Big[1+\frac n2\left(1-\frac 1q\right)+\frac 1{w^2}\Big(1-\frac n2\Big(1-\frac1q\Big)\Big)\Big]^q\omega^{2q-n(q-1)}-1
\end{equation}
there exists a ground state $(u^*,v^*)$ with $u^*,v^*>0$.

\bigskip

In the present paper we will prove the existence of a positive radial decreasing ground state solution to (\ref{elliptic}) for all $(q,n)$ satisfying the condition (\ref{cond}). Exploring this radial decay, we  improve the constant $\mathcal{C}_{\omega,n,q}$ derived in \cite{MMP} for all $q>1$ and large $\omega$ in the case $n=1$ and for all $1<q<2$ in any dimension, in fact replacing it by $0$ in the latter case.

\bigskip

When dealing with the system (\ref{elliptic}) it is often necessary to treat the case $n=1$ separately due to the lack of compactness of the injection $H^1_d(\er)\hookrightarrow L^q(\er)$, $q>2$, where $H_d^1(\er)$ denotes the space of the radially symmetric functions of  $H^1(\er)$. This lack of compactness is, in a sense, a consequence of the inequality
\begin{equation}
 \label{ineq}
|u(x)|\leq C|x|^{\frac{1-n}2}\|u\|_{H^1(\er^n)}
\end{equation}
for $u\in H_d^1(\er)$. Indeed, (\ref{ineq}) gives no decay in the case $n=1$. However, if $u$ is also radially decreasing, it is easy to establish that
$$|u(x)|\leq C|x|^{-\frac{n}2}\|u\|_{L^2(\er^n)},$$
which provides decay in all space dimensions, hence compacity by applying the classical Strauss' compactness lemma (\cite{Strauss}). Hence, putting $$H_{rd}^1(\er^n)=\{u \in H^1_d(\er^n)\,:\,u\textrm{ is radially decreasing}\},$$ we get the compactness of the injection $H_{rd}^1(\er^n)\hookrightarrow L^q(\er^n)$ for all $n\geq 1$ (see the Appendix of \cite{BL} for more details). We will use this fact to present a unified approach for the problem of the energy minimization of (\ref{elliptic}), valid in all space dimensions.

\bigskip
Before stating our results more precisely, and following the functional settings in \cite{Colorado}, \cite{Coloradobis} and \cite{MMP}, let us introduce a few notations: we denote by $\|\cdot\|_q$ the standard $L^q(\er^n)$ norm and, for $(u,v)\in E:=H^1(\er^n)\times H^1(\er^n)$, we put
$$\|(u,v)\|_{\dot\omega}^2:=\|u\|^2+\|v\|_{\dot\omega}^2:=\|u\|_2^2+\|\nabla u\|_2^2+\omega^2\|v\|_2^2+\|\nabla v\|_2^2.$$
We introduce the energy functional associated to (\ref{elliptic}),
$$I(u,v):=\frac 12\|(u,v)\|_{\dot \omega}^2-\frac 1{2q}\Big(\|u\|_{2q}^{2q}+\|v\|_{2q}^{2q}+2b\|uv\|_q^q\Big),$$
noticing that $(u,v)$ is a solution of (\ref{elliptic}) if and only if $\nabla I(u,v)=0$.\\ We will study the minimization problem
\begin{equation}
 \label{minimization}
\inf \{I(u,v)\,:\,(u,v)\in \mathcal{N}\},
\end{equation}
where the so-called Nehari manifold $\mathcal{N}$ is defined by
$$\mathcal{N}:=\{(u,v)\in H^1(\er^n)\times H^1(\er^n)\,: (u,v)\neq (0,0), \nabla I(u,v)\perp (u,v)\},$$
that is, $(u,v)\in \mathcal N$ if and only if $(u,v)\neq (0,0)$ and
$$\tau(u,v):=\langle \nabla I(u,v),(u,v)\rangle_{L^2}=\|(u,v)\|_{\dot \omega}^2-\Big(\|u\|_{2q}^{2q}+\|v\|_{2q}^{2q}+2b\|uv\|_q^q\Big)=0.$$
As pointed out in \cite{Colorado} for the case $q=2$, we notice that 
$$\langle \nabla \tau(u,v),(u,v) \rangle_{L^2}=2\|(u,v)\|_{\dot \omega}^2-2q\Big(\|u\|_{2q}^{2q}+\|v\|_{2q}^{2q}+2b\|uv\|_q^q\Big),$$
and, if $(u,v)\in\mathcal{N}$, 
\begin{equation}
\label{estmanifold}
\langle \nabla \tau(u,v),(u,v) \rangle_{L^2}=2(1-q)\|(u,v)\|_{\dot \omega}^2<0
\end{equation}
which shows that $\mathcal{N}$ is locally smooth.\\ Furthermore, it is easy to check that $[h_1,h_2]\textrm{Hess}\,\tau_{(0,0)}\,^t[h_1,h_2]>0$ for all $(h_1,h_2)\neq(0,0)$: $(0,0)$ is a strict minimizer of $\tau$, hence an isolated point of the set $\{\tau(u,v)=0\}$, implying that $\mathcal{N}$ is a complete manifold. Finally, any critical point of $I$ constrained to $\mathcal{N}$ is a critical point of $I$. Indeed, let us consider $(u,v)\in\mathcal{N}$ a critical point of $I$ constrained to $\mathcal{N}$. There exists a Lagrange multiplier $\lambda$ such that $\nabla I(u,v)=\lambda \nabla \tau(u,v).$\\
By taking the $L^2$ scalar product with $(u,v)$, $$\langle \nabla I(u,v),(u,v)\rangle_{L^2}=\lambda \langle\nabla \tau(u,v),(u,v)\rangle_{L^2},$$ that is, in view of (\ref{estmanifold}), $0=\lambda (2-2q)\|(u,v)\|_{\dot \omega}^2$, hence $\lambda=0$ and $\nabla I(u,v)=0$.

\bigskip

Putting $E_{rd}=H_{rd}^1\times H_{rd}^1$ the cone of symmetric radially decreasing non-negative  functions of $E$, we will prove the following result:
\begin{Teorema}
\label{Teorema1}
Let $n\geq 1$ and $q>1$, with $q\leq \frac{n}{n-2}$ if $n> 3$. There exists a minimizing sequence $(u_n,v_n)\in E_{rd}$ for the minimization problem (\ref{minimization}). Furthermore, 
$(u_n,v_n)\to (u_*,v_*)\in E_{rd}$ strongly in $H^1(\er^n)\times H^1(\er^n)$. In particular
\begin{multline}
I(u_*,v_*)=\min_{\mathcal N} I(u,v)=\min_{\mathcal{N}\cap E_{rd}}I(u,v)\\=\min \{I(u,v)\,:\,(u,v)\neq(0,0)\textrm{ and }\nabla I(u,v)=0\}.
\end{multline}
\end{Teorema}
Concerning the existence of ground states with non-trivial components, we will show:
\begin{Teorema}
 \label{Teorema2}
Let $n\geq 1$ and $1<q<2$, with $q<\frac{n}{n-2}$ if $n\geq 3$.\\
Then for all $b>0$ there exists a ground state solution $(u,v)\in E_{rd}$ to (\ref{elliptic}) with $u>0$ and $v>0$.
\end{Teorema}
\begin{Teorema}
 \label{Teorema3}
Let $n=1$ and $q\geq 2$. If
\begin{equation}
 \label{conditionb}
 b\geq \mathcal{D}_{\omega,q}=\frac{2^q-1}2\omega^{1+\frac q2}-\frac12\omega^{-\frac q2}
\end{equation}
there exists a ground state solution $(u,v)\in E_{rd}$ to (\ref{elliptic}) with $u>0$ and $v>0$.
\end{Teorema}
Notice that
$$\mathcal{D}_{\omega,q}<\mathcal{C}_{\omega,1,q}=\frac 12\Big(\frac 32-\frac 1{2q}+\frac 1{\omega^2}\Big(\frac 12+\frac1 {2q}\Big)\Big)^q\omega^{1+q}-1$$
at least for large values of $\omega$.
\section{Proof of Teorem \ref{Teorema1}}
We begin by observing that for $(f,g)\in E$, $(f,g)\neq (0,0)$, with $\tau(f,g)\leq 0$, there exists $t\in]0,1]$ such that $(tf,tg)\in\mathcal{N}$. Indeed, if $\tau(f,g)=0$, we choose $t=1$. If $\tau(f,g)<0$ we simply notice that $$\tau(tf,tg)=t^2\Big(\|(f,g)\|_{\dot \omega}^2-t^{2q-2}(\|f\|_{2q}^{2q}+\|g\|_{2q}^{2q}+2b\|fg\|_q^q)\Big):=t^2T_{f,g}(t),$$
with $T_{f,g}(0)>0$ and $T_{f,g}(1)<0$.\\
Also, we notice that if $(f,g)\in\mathcal{N}$,
\begin{equation}
 \label{valori}
I(f,g)=\Big(\frac 12-\frac 1{2q}\Big)\|(f,g)\|_{\dot \omega}=\Big(\frac 12-\frac 1{2q}\Big)(\|f\|_{2q}^{2q}+\|g\|_{2q}^{2q}+2b\|fg\|_q^q).
\end{equation}
We now take a minimizing sequence $(u_n,v_n)\in \mathcal{N}$ for the problem
$$m=\inf\{I(u,v)\,:\,(u,v)\in\mathcal{N}\}.$$
From (\ref{valori}), it is clear that $m\geq 0$ and that  $(u_n,v_n)$ is bounded in $E$.

\bigskip

We put $u_n^*$ and $v_n^*$ the decreasing radial rearrangements of $|u_n|$ and $|v_n|$ respectively. It is well-known that this rearrangement preserves the $L^p$ norm ($1\leq p \leq +\infty$). Furthermore, the P\'olya-Szeg\"o inequality
$$\|\nabla f^*\|_2\leq \|\nabla |f|\|_2$$ in addition with the inequality $\|\nabla |f|\|_2\leq \|\nabla f\|_2$ (see \cite{Lions1}) shows that
$$\|(u_n^*,v_n^*)\|_{\dot \omega}^2\leq \|(u_n,v_n)\|_{\dot \omega}^2.$$
On the other hand, the Hardy-Littlewood inequality
$$\int |fg|\leq \int f^*g^*$$
combined with the monotonicity of the map $\lambda\to\lambda^q$ (see for instance \cite{Rear} for details) yields $\|fg\|_q\leq \|f^*g^*\|_q$ and, finally,
$$\tau(u_n^*,v_n^*)\leq \tau(u_n,v_n)=0.$$
Next, let $t_n\in]0,1]$ such that $(t_nu_n^*,t_nv_n^*)\in \mathcal{N}.$
We obtain 
$$I(t_nu_n^*,t_nv_n^*)=t_n^2\Big(\frac 12-\frac 1{2q}\Big)\|(u_n^*,v_n^*)\|_{\dot \omega}^2\leq \Big(\frac 12-\frac 1{2q}\Big)\|(u_n,v_n)\|_{\dot \omega}^2=I(u_n,v_n)$$ and we obtained a minimizing sequence $(t_nu_n^*,t_nv_n^*)$ in $E_{rd}$, denoted again, in what follows, by $(u_n,v_n)$.
Since this sequence is bounded in $H^1(\er^n)$, up to a subsequence, $(u_n,v_n)\rightharpoonup (u_*,v_*)$ in $H^1(\er^n)$ weak. Also, since the injection $E_{rd}\to L^{2q}(\er^n)$ is compact, up to a subsequence, $(u_n,v_n)\to (u_*,v_*)$ in $L^{2q}(\er^n)$ strong.\\
Hence, since $\|u_n\|_{2q}^{2q}+\|v_n\|_{2q}^{2q}+2b\|u_nv_n\|_q^q\to \|u_*\|_{2q}^{2q}+\|v_*\|_{2q}^{2q}+2b\|u_*v_*\|_q^q$, we deduce that
$$\tau(u_*,v_*)\leq \liminf \tau(u_n,v_n)=0.$$
Once again, let $t\in]0,1]$ such that $(tu_*,tv_*)\in \mathcal{N}$.
\begin{multline*}
m\leq I(tu_*,tv_*)=t^2\Big(\frac 12-\frac 1{2q}\Big)\|(u_*,v_*)\|_{\dot \omega}^2\\
\leq \Big(\frac 12-\frac 1{2q}\Big)\liminf \|(u_n,v_n)\|_{\dot \omega}^2\leq \liminf I(u_n,v_n)=m.
\end{multline*}
This implies that $(tu_*,tv_*)$ is a minimizer. In particular, all inequalities above are in fact equalities: $t=1$, $(u_*,v_*)\in \mathcal{N}$, $\|(u_*,v_*)\|_{\dot \omega}=\lim \|(u_n,v_n)\|_{\dot \omega}$, 
$\|u_n\|_{H^1}\to \|u_*\|_{H^1}$, $\|v_n\|_{H^1}\to \|v_*\|_{H^1}$ and $(u_n,v_n)\to (u_*,v_*)$ in $H^1(\er^n)$ strong.\\
Finally, it is clear that $(u_*,v_*)$ is a ground state:
if $(w_1,w_2)\neq(0,0)$ is a critical point of $I$ such that $I(w_1,w_2)<I(u_*,v_*)$, taking once again $w_1^*$ and $w_2^*$ the decreasing radial rearrangements of $|w_1|$ and $|w_2|$, there exists $t\in]0,1]$ such that $(tw_1^*,tw_2^*)\in\mathcal{N}$ and $I(tw_1^*,tw_2^*)\leq I(w_1,w_2)$, which leads to a contradiction. This completes the proof of Theorem \ref{Teorema1}.  \hfill$\blacksquare$

\section{Ground states with non-trivial components}
Let $(u_*,v_*)\in E_{rd}$ the ground state mentionned in Theorem \ref{Teorema1}. If $v_*=0$, $u_*=u_0$ is the unique positive radially symmetric solution of the elliptic equation $-\Delta u+u=u^{2q-1}$ (see \cite{Unic}).\\
Also, if $u_*=0$, $v_*=v_0$ is the unique positive radially symmetric solution of $-\Delta v+\omega^2v=v^{2q-1}$, which relates to $u_0$ by the relation $v_0(x)=\omega^{\frac 1{q-1}}u_0(\omega x)$.\\
Hence, to show the existence of a ground state with nontrivial components, we only have to exhibit an element $(f,g)\in \mathcal{N}\cap E_{rd}$, $f\neq 0$, $g\neq 0$, such that
\begin{equation}
\label{eninf}
 I(f,g)\leq \min\{I(u_0,0),I(0,v_0)\}.
\end{equation}
\noindent
Since $I(u_0,0)=\Big(\frac 12-\frac 1{2q}\Big)\|u_0\|_{2q}^{2q}$, $I(0,v_0)=\omega^{\frac{2q}{q-1}-n}\Big(\frac 12-\frac 1{2q}\Big)\|u_0\|_{2q}^{2q}$ and\\ $\frac{2q}{q-1}-n>0$, for $\omega \geq 1$ the inequality (\ref{eninf}) reduces to
\begin{equation}
 \label{eninf2}
 I(f,g)\leq I(u_0,0).
\end{equation}

\bigskip

\noindent
We first compute $x>0$ such that $(f,g):=(x u_0,x\theta v_0)\in \mathcal{N}$, where $\theta>0$ will be chosen later (see \cite{Col1} and \cite{Col} for a recent application of a related technique to the Schr\"odinger-KdV system):
$$\tau(f,g)=x^2\|(u_0,\theta v_0)\|_{\dot \omega}^2-x^{2q}\Big(\|u_0\|_{2q}^{2q}+\theta^{2q}\|v_0\|_{2q}^{2q}+2b\theta^q\|u_0v_0\|_q^q\Big)=0.$$
Since 
$$\|\theta v_0\|_{\dot \omega}^2=\omega^{2+\frac 2{q-1}-n}\|\theta u_0\|_2^2+\omega^{2+\frac 2{q-1}-n}\|\theta \nabla u_0\|_2^2=\omega^{\frac {2q}{q-1}-n}\theta^2 \|u_0\|^2$$
and 
$$\|v_0\|_{2q}^{2q}=\omega^{\frac{2q}{q-1}-n}\|u_0\|_{2q}^{2q},$$
we obtain
$$x^{2q-2}=\frac{(1+\theta^2\omega^{\frac {2q}{q-1}-n})\|u_0\|^2}{(1+\theta^{2q}\omega^{\frac {2q}{q-1}-n})\|u_0\|_{2q}^{2q}+2b\theta^q\|u_0v_0\|_q^q}=\frac{1+\theta^{2}\omega^{\frac {2q}{q-1}-n}}{1+\theta^{2q}\omega^{\frac {2q}{q-1}-n}+2b\theta^q\frac{\|u_0v_0\|_q^q}{\|u_0\|_{2q}^{2q}}}.$$
Since $u_0$ is radial and nonincreasing and $\omega\geq 1$, 
$$\|u_0v_0\|_q^q=\omega^{\frac q{q-1}}\int u_0^q(x)u_0^q(\omega x)dx\leq \omega^{\frac q{q-1}}\int u_0^q(x)u_0^q(x)dx=\omega^{\frac q{q-1}}\|u_0\|_{2q}^{2q}.$$
Also,
$$\|u_0v_0\|_q^q\geq \omega^{\frac q{q-1}}\int u_0^{2q}(\omega x)dx= \omega^{\frac q{q-1}-n}\|u_0\|_{2q}^{2q}.$$
Hence, we obtain
\begin{equation}
\label{x}
\frac{1+\theta^2\omega^{\frac {2q}{q-1}-n}}{1+\theta^{2q}\omega^{\frac {2q}{q-1}-n}+2b\theta^q\omega^{\frac q{q-1}}}\leq x^{2q-2}\leq\frac{1+\theta^2\omega^{\frac {2q}{q-1}-n}}{1+\theta^{2q}\omega^{\frac {2q}{q-1}-n}+2b\theta^q\omega^{\frac q{q-1}-n}}
\end{equation}
and
$$I(f,g)=x^2\Big(\frac 12-\frac 1{2q}\Big)\|(u_0,\theta v_0)\|_{\dot \omega}^2=x^{2}\Big(\frac 12-\frac 1{2q}\Big)(1+\theta^2\omega^{\frac {2q}{q-1}-n})\|u_0\|^2.$$
The condition (\ref{eninf2}) then becomes $x^2(1+\theta^2\omega^{\frac {2q}{q-1}-n})\leq 1.$\\
In view of (\ref{x}), a sufficient condition is
$$\frac{(1+\theta^2\omega^{\frac {2q}{q-1}-n})^q}{1+\theta^{2q}\omega^{\frac {2q}{q-1}-n}+2b\theta^q\omega^{\frac q{q-1}-n}}\leq 1,$$
that is,
$$b\geq \frac{(1+\theta^2\omega^{\frac {2q}{q-1}-n})^q-1-\theta^{2q}\omega^{\frac {2q}{q-1}-n}}{2\theta^q\omega^{\frac q{q-1}-n}}.$$
We now put $\theta^2=\epsilon^2\omega^{n-\frac {2q}{q-1}}$, for $\epsilon>0$, obtaining the condition
$$b\geq \frac{(1+\epsilon^2)^q-1}{2\epsilon^q}\omega^{q-\frac n2(q-2)}-\frac 12\epsilon^q\omega^{(\frac n2-1)q}.$$
For $1<q<2$, $\displaystyle \lim_{\epsilon\to 0} \frac{(1+\epsilon^2)^q-1}{2\epsilon^q}=0$.\\
Hence, the arbitrary value of $\epsilon$ establishes the sufficient condition $b>0$.

\bigskip

\noindent
For $n=1$, putting $\epsilon=1$, we obtain the bound 
\begin{equation}
 \label{melhorar}
b\geq \frac{2^q-1}2\omega^{1+\frac q2}-\frac12\omega^{-\frac q2},
\end{equation}

as  stated in Theorem \ref{Teorema3}.\hfill$\blacksquare$\\
We finish by making a few remarks:
\begin{Rem}
 For $\omega=1$ and $\theta=1$, we obtain, for all $n\geq 1$, the bound $2^{q-1}-1$ which is known to be optimal for $q\geq 2$, in the sense that for $b<2^{q-1}-1$ all ground states of (\ref{elliptic}) have one null component  (see \cite{MMP}, Theorem 2.5).
\end{Rem}
\begin{Rem}
The bound in (\ref{melhorar}) can be slightly improved for large values of $\omega$ by replacing the quantity $\displaystyle \frac{2^q-1}2$ by the minimum of $\displaystyle \frac{(1+\epsilon^2)^q-1}{2\epsilon^q}$ for $\epsilon>0$.
\end{Rem}
\begin{Rem}
 For $n\geq 4$ we have $1<q<2$, hence the problem of the existence of ground states with non-trivial components is completely solved for these spatial dimensions.
\end{Rem}
\noindent
{\bf Acknowledgment:} This article was partially supported by Funda\c{c}\~ao para a Ci\^encia e Tecnologia, through contract PEst-OE/MAT/UI0297/2015.

\end{document}